\documentclass[11pt,reqno]{amsart}

\title[]
{Magnetically Forced Electroconvection on Bounded Domains}

\author[E. Abdo]{Elie Abdo}
\address[E. Abdo]
{
    Department of Mathematics \\
    American University of Beirut \\
    Beirut 1107-2020\\
    Lebanon.
} 
\email{ea94@aub.edu.lb}

\author[J. Germany]{Joe Germany}
\address[J. Germany]
{
    Department of Mathematics \\
    American University of Beirut \\
    Beirut 1107-2020\\
    Lebanon.
}
\email{jmg15@mail.aub.edu}

\author[Q. Lin]{Quyuan Lin}
\address[Q. Lin]
{
    School of Mathematical and Statistical Sciences \\
	Clemson University \\
	Clemson, SC 29634, USA.
} 
\email{quyuanl@clemson.edu}

\usepackage{xfrac}

\usepackage[margin=1in]{geometry}
\usepackage{amsmath, amsthm, amssymb, yhmath, xfrac, mathtools}
\usepackage[most]{tcolorbox}
\usepackage[normalem]{ulem}

\usepackage{times}
\usepackage{color}
\usepackage{hyperref}
\newcommand{\pa}{\partial}
\newcommand{\la}{\label}
\newcommand{\fr}{\frac}
\newcommand{\na}{\nabla}
\newcommand{\be}{\begin{equation}}
\newcommand{\ee}{\end{equation}}
\newcommand{\ba}{\begin{array}{l}}
\newcommand{\ea}{\end{array}}

\newtheorem{thm}{Theorem}[section]
\newtheorem{prop}[thm]{Proposition}
\newtheorem{lemma}[thm]{Lemma}

\newcommand{\G}{\Gamma}

\newcommand{\eps}{\epsilon}

\newcommand{\D}{\Delta}
\newcommand{\Div}[1]{\na \cdot #1}

\newcommand{\OInt}[1]{\int_{\Omega} #1 \, dx}
\newcommand{\ComInt}[1]{\int_0^\infty t^{-1+\fr{s}{4}} \int_{\Omega} #1 \, dy \, dt}

\newcommand{\TInt}[1]{\int_0^t #1 \, ds}

\newcommand{\dt}[0]{\fr{d}{dt}}
\newcommand{\Lnorm}[2]{\left\lVert #1 \right\rVert_{L^{#2}}}
\newcommand{\LnormD}[3]{\left\lVert #1 \right\rVert_{L^{#2} \left({#3}\right)}}

\renewcommand{\l}{\Lambda}

\usepackage{cite}

\usepackage{MnSymbol,wasysym}
\usepackage{accents}

\newcommand{\N}{\mathbb N}

\newcommand{\R}{\mathbb R}

\def\PP{\mathbb P}

\date{\today}
\begin{document}
\begin{abstract} 
We study a magnetically forced electroconvection model describing a charged smectic thin film coupled to three-dimensional electrostatics. We establish global existence, uniqueness, regularity, and exponential decay of solutions on bounded domains for arbitrary magnetic field strength, assuming initial velocity in $H^1$ and initial charge density in $L^4$. The analysis relies on a new hierarchy of coupled energies and new commutator estimates for negative powers of fractional Laplacians on bounded domains. These commutator estimates provide a general framework for studying low-regularity well-posedness questions for nonlinear PDEs on bounded domains.
\end{abstract} 

\noindent

\maketitle

\section{Introduction}

We consider a thin smectic film occupying a two-dimensional (2D) domain $\Omega \subset \R^2$, embedded in the three-dimensional (3D) space $\Omega \times \R$. The charge carriers are confined to the film, while the surrounding medium remains electrically active. This dimensional contrast of charges evolving on a 2D surface while generating an electric field in 3D gives rise to an electroconvection model with nonlocal aspects. 

Let $q(x,t)$ denote the surface charge density and $u(x,t) = (u_1(x,t), u_2(x,t))$ the velocity of the incompressible fluid in the film. The charge density satisfies the conservation law
\[
\partial_t q + \nabla \cdot j = 0,
\]
where the surface current density $j$ is obtained by restricting a 3D current density $\tilde j$ to the plane. In the presence of a perpendicular magnetic field $\tilde B=(0,0,B)$, Ohm's law together with the Lorentz correction give
\[
\tilde j = (u_1,u_2,0)q + \sigma\big(\tilde E + (u_1,u_2,0)\times \tilde B\big),
\]
where $\sigma>0$ is the conductivity of the medium and $\tilde{E} = (E_1, E_2, E_3)$ is the 3D electric field generated by the charge density. Restricting to $\Omega$ yields
\[
j = uq + \sigma\big(E - B u^\perp\big),
\]
where $u^\perp = (-u_2,u_1)$ and $E= (E_1, E_2)$.
The electric field $\tilde E$ is determined in the full 3D space by Gauss' law
\[
-\eps_0 \nabla_3 \cdot \tilde E = 2q\,\delta_\Omega,
\]
which reflects that the charge density is supported only on the film. {Here $\eps_0>0$ is the vacuum permittivity, $\delta_{\Omega}$ is the Dirac mass on $\Omega$, and $\na_3$ is the 3D gradient operator.} Solving this equation and restricting back to $\Omega$ leads to the nonlocal relation
\[
E = -\frac{1}{\eps_0}\nabla \Lambda^{-1} q,
\]
where $\Lambda = (-\Delta)^{1/2}$ with homogeneous Dirichlet boundary conditions. This operator arises precisely from the mismatch between the 3D electrostatics and the 2D dynamics.
Substituting the expression for $j$ into the conservation law and using the incompressibility of the fluid, we obtain
\[
\partial_t q + u \cdot \nabla q + \frac{\sigma}{\eps_0}\Lambda q - \sigma B \nabla \cdot u^\perp = 0.
\]
The term $\Lambda q$ represents a nonlocal dissipative mechanism induced by the electric field and the surrounding medium.

The fluid motion is governed by the incompressible Navier--Stokes equations with viscosity $\nu>0$. The forces acting on the fluid arise from both electrostatic and magnetic effects. The electric force is given by $qE$, while the magnetic contribution is described by the Lorentz force
\[
(\tilde j \times \tilde B)\big|_\Omega = -B j^\perp = -B (uq)^\perp - \sigma B E^\perp - \sigma B^2 u,
\]
yielding the total forcing
\[
qE + (\tilde j \times \tilde B)\big|_\Omega
= -\frac{1}{\eps_0} q R q - B q u^\perp + \frac{\sigma B}{\eps_0} R^\perp q - \sigma B^2 u,
\] where $R = \na \Lambda^{-1}$ is the 2D Riesz transform operator. Consequently, the momentum equation is described by
\[
\partial_t u + u \cdot \nabla u + \nabla p - \nu \Delta u
= -\frac{1}{\eps_0} q R q - B q u^\perp + \frac{\sigma B}{\eps_0} R^\perp q - \sigma B^2 u.
\] Combining the above relations gives rise to the coupled nonlinear and nonlocal system 
\begin{subequations}\label{main-system}
\begin{align}
    &\partial_t q + u \cdot \nabla q + \frac{\sigma}{\eps_0}\Lambda q - \sigma B \nabla \cdot u^\perp = 0, \label{charge_density_eq} \\
    &\partial_t u + u \cdot \nabla u + \nabla p - \nu \Delta u
    = -\frac{1}{\eps_0} q R q - B q u^\perp + \frac{\sigma B}{\eps_0} R^\perp q - \sigma B^2 u, \label{navier_stokes_eq} \\
    &\nabla \cdot u = 0 \label{div_free_eq},
\end{align}
\end{subequations}
supplemented with homogeneous Dirichlet boundary conditions $(u,q)|_{\partial\Omega}=0$. 
This system reflects the interplay between transport, nonlocal dissipation induced by the 3D electrostatics, and the combined electrostatic and magnetic forces acting on the fluid confined to the thin film. This paper is devoted to the global well-posedness and long-time dynamics of the aforementioned electroconvection model. For simplicity, we take $\sigma=\epsilon_0=1$ because these parameters do not affect the analysis of the problem. However, we keep track of the magnetic component $B$ as we are interested in its effect on the long-time behavior of solutions.

Electroconvection models without magnetic effects, derived from experiments in annular smectic thin films \cite{constantin2017some, daya1999electrically, daya1998annular, tsai2007direct}, were shown to be globally well-posed on bounded domains in high-regularity Sobolev spaces \cite{constantin2017some}. Under periodic boundary conditions, global well-posedness, uniqueness, and the existence of finite-dimensional global attractors were established using nonlinear maximum principles and commutator estimates \cite{abdo2021long}. In the whole-space setting $\mathbb{R}^2$, sharp decay rates for the charge density and velocity were obtained and justified through comparisons with linearized systems and Fourier-based methods \cite{abdo2024long}. On bounded domains with Dirichlet boundary conditions, global well-posedness for $H^1$ initial data, instantaneous smoothing, and exponential decay were proved via spectral regularization of nonlocal operators, along with the existence of smooth global attractors in forced settings \cite{abdo2023long}. On the other hand, it has been shown in experimental studies that applied transversal magnetic fields significantly impact the dynamics of electroconvection phenomena \cite{gleeson1996onset, ahmed2023effect, huh1999formation, huh2007electrohydrodynamic, huh2015multiplicative}. In particular, when a magnetic field is applied,  the critical voltage required for instability increases, leading to enhanced system stability \cite{ahmed2023effect}.
Very recently, electroconvection in porous media with magnetic effects was analyzed, yielding global weak solutions for large $L^\infty$ data and convergence toward quasigeostrophic regimes under suitable scalings \cite{abdo2025electroconvection}.

In this paper, we address the magnetically forced electroconvection model coupled to the Navier-Stokes flow, and we construct unique global solutions that decay exponentially in time to zero under the assumptions of $H^1$ initial velocity, $L^4$ initial charge density, and an external magnetic field of arbitrary size. The argument proceeds in three main stages. In the first stage, we couple the $H^{-1/2}$ evolution of the charge density with the $L^2$ evolution of the velocity. This coupling yields a crucial cancellation of the nonlinear contribution arising from the electrostatic forcing term $q R q$. In the second stage, we exploit this structure to propagate decay, upgrading the control of the charge density from $H^{-1/2}$ to $L^2$. This step also provides cancellation of the nonlinear advection term $u \cdot \nabla q$ and relies fundamentally on the dissipative mechanism identified in the first stage. Finally, in the third stage, we couple the $H^1$ evolution of the velocity with the $L^4$ evolution of the charge density. This step builds on the dissipative structures established in the previous stages and allows us to close the argument. In contrast with previously studied electroconvection models without magnetic effects, where no such coupling was necessary, the presence of the magnetic field enforces a nontrivial coupling between the energies, both in the construction of solutions and in the analysis of their long-time behavior. 

The uniqueness of solutions in functional spaces weaker than \( H^1 \) is particularly challenging due to the nonlinear transport structure, which inherently involves gradients through divergence-type terms. In this paper, we establish global well-posedness for initial charge density in \( L^4 \) and initial velocity in \( H^1 \), both in the presence and in the absence of an externally imposed magnetic field. This improves upon previous results obtained on bounded domains that required \( H^1 \) regularity for both the velocity and the charge density. To achieve this, we introduce a suitable topology in which the difference between two solutions vanishes whenever their initial data coincide. This framework relies on negative fractional Sobolev spaces, leading to a delicate analysis of the nonlinear advection term \( u \cdot \nabla q \). Exploiting the divergence-free condition on the velocity field, we construct a commutator that cancels the highest-order contributions arising from the nonlinearity, and we derive precise estimates to control it.

A key step in our approach is addressing a more general problem: establishing commutator estimates for negative powers of the Laplacian on bounded domains. {Compared to the whole spaces or periodic tori, this task is nontrivial on bounded domains} due to the lack of Fourier-based tools, the nonlocal nature of fractional Laplacians, and their failure to commute with differential operators. Our analysis overcomes these difficulties and yields commutator bounds that depend explicitly on the distance-to-the-boundary function. By invoking Rellich-type inequalities, we obtain estimates that remain uniform up to the boundary, rather than blowing up near it. Such estimates constitute a main novelty of the paper and provide a robust {general} framework for studying well-posedness of PDEs in low-regularity Sobolev spaces.

This paper is organized as follows. In Section~\ref{sec:pre}, we introduce some preliminaries and notation that will be used throughout the paper. In Section~\ref{sec_commutator_estimate}, we establish new commutator estimates for negative powers of the fractional Laplacian on bounded domains. In Section~\ref{sec_existence_and_uniqueness}, we address the existence, uniqueness, and long-time behavior of solutions to the magnetically forced electroconvection model based on the aforementioned commutator estimates. Appendix~\ref{kernel_appendix} is devoted to Dirichlet heat kernel estimates that are required for the analysis of commutators. Finally, Appendix~\ref{gronwall_section} is dedicated to a Gr\"onwall-type lemma that is employed to derive the exponential decay in time of solutions.

\section{Preliminaries}\label{sec:pre}

 Let $\Omega\subset \mathbb R^2$ be a bounded domain with a smooth boundary.
 
{\bf{Fractional Powers of the Laplacian}}. We denote by $\Delta$ the Laplacian operator with homogeneous Dirichlet boundary conditions. We note that $-\Delta$  is defined on $\mathcal{D}(-\Delta) = H^2(\Omega) \cap H_0^1(\Omega)$, and is positive and self-adjoint in $L^2(\Omega)$. We consider an orthonormal basis of $L^2(\Omega)$ consisting of eigenfunctions $\left\{w_j\right\}_{j=1}^{\infty} \subset H_0^1(\Omega)$ of $-\Delta$ satisfying
\be
-\Delta w_j = \lambda_j w_j,
\ee
where the eigenvalues $\lambda_j$ obey $0 < \lambda_1 \leq ... \leq \lambda_j \le ...\rightarrow \infty$. 
For $s \in \R$, we define the fractional Laplacian operator of order $s$, denoted by $\l^s$, by
\be \la{maindef}
\l^s f= \sum_{j=1}^{\infty} \lambda_j^{\fr{s}{2}} (f, w_j)_{L^2} w_j,
\ee with domain 
\be 
\mathcal{D}(\l^s) = \left\{f \in L^2(\Omega) : \|\l^{s} f\|_{L^2}^2 := \sum\limits_{j \in \N} \lambda_j^{s}(f, w_j)_{L^2}^2 < \infty \right\}.
\ee 
For $s \in [0,1]$, we identify the domains $\mathcal{D}(\l^s)$ with the usual Sobolev spaces as follows, 
\be 
\mathcal{D}(\l^s) = \begin{cases} H^{s}(\Omega) \hspace{8cm} \mathrm{if} \; s \in [0, \fr{1}{2}) 
\\ H_{00}^{\fr{1}{2}} (\Omega) = \left\{f \in H_0^\fr{1}{2} (\Omega) : f/\sqrt{d(x)} \in L^2(\Omega)\right\} \hspace{2cm} \mathrm{if} \; s = \fr{1}{2}
\\H_{0}^{s} (\Omega) \hspace{8cm} \mathrm{if} \; s \in (\fr{1}{2},1] 
\end{cases}
\ee where 
$d(x)$ is the distance-to-the-boundary function. 

{\bf{Stokes Operator.}} We denote by $H$ and $V$ the spaces
\be 
H = \left\{v \in (L^2(\Omega))^2: \na \cdot v = 0, v \cdot n|_{\pa \Omega} = 0 \right\},
\ee where $n$ is the outward unit normal to $\pa \Omega$, and 
\be 
V = \left\{v \in (H_0^1(\Omega))^2 : \na \cdot v = 0\right\}.
\ee 
Let $\PP: (L^2(\Omega))^2 \rightarrow H$ be the Leray Hodge projection. We define the Stokes operator, denoted by $A$, as 
\be 
A:= - \PP \Delta,
\ee with domain $\mathcal{D}(A) = V \cap (H^2(\Omega))^2$. $A$ is positive, self-adjoint, and injective, and its inverse $A^{-1}$ is compact. We denote the eigenvalues of $A$ by $\mu_j, j = 1,2,...$, and the corresponding eigenfunctions by $\phi_j, j = 1,2...$, and we note that $0 < \mu_1 \le ... \le \mu_j \le ... \rightarrow \infty$. We define the fractional powers of the Stokes operator, denoted by $A^s$, as
\be \la{stoeig}
A^s v= \sum_{j=1}^{\infty} \mu_j^s (v, \phi_j)_{L^2} \phi_j,
\ee with domain 
\be 
\mathcal{D}(A^s) = \left\{v \in H : \|A^{s} v\|_{L^2}^2 := \sum\limits_{j \in \N} \mu_j^{2s}(v, \phi_j)_{L^2}^2 < \infty \right\}.
\ee 

{\bf{Notation.}} {Throughout the paper, the universal constants $c$ and $C$ may change from line to line. When necessary, we use subscripts to indicate the dependence of the constant on certain parameters, {\it e.g.}, we write $C_r$ to emphasize that the constant depends on $r$.}

\section{Commutator Estimates for the Negative Fractional Laplacian}
\label{sec_commutator_estimate}

The representation of negative powers of the Laplacian using the heat kernel is given by
\begin{equation}
\label{negative_representation}
    \l^{-s} f(x) = c_s \int_0^\infty t^{-1+\fr{s}{2}} \,e^{t\D} f(x) \ dt,
\end{equation}
for $s > 0$, where
\begin{equation}
    e^{t \D} f(x) = \int_\Omega H(x,y,t) f(y) \ dy, \nonumber
\end{equation}
for all $x \in \Omega$.
From \cite{li1986parabolic}, we know that the following bounds on the gradient of the Dirichlet heat kernel hold
\be\begin{aligned}
\label{gradient_bound}
    |\na_x H_D| &\leq C t^{-\fr{1}{2} - \fr{d}{2}} e^{-\fr{|x-y|^2}{Kt}}, \\
    |\na_y H_D| &\leq C t^{-\fr{1}{2} - \fr{d}{2}} e^{-\fr{|x-y|^2}{Kt}},
\end{aligned}\ee
for all $(x,y) \in \Omega \times \Omega$ and $t > 0$. \\
The following standard estimate will be frequently used in this section:
\begin{equation}
    \label{standard_integral_estimate}
    \int_0^\infty t^{- 1 - \fr{m}{2}} e^{-\fr{p^2}{Kt} }\ dt \leq C_{K, m} \, p^{-m},  \qquad m, \,p, \, K >0.
\end{equation}
For $0 \leq \alpha \leq 1$, we recall the H\"older space $C^{0, \alpha} (\overline{\Omega})$ with the associated norm
\begin{equation}
\label{holder_norm}
    \lVert f \rVert_{C^{0, \alpha}} = \Lnorm{f}{\infty} + [f]_{C^{0, \alpha}},
\end{equation}
where the seminorm is defined as
\begin{equation}
\label{holder_seminorm}
    [f]_{C^{0,\alpha}} = \sup_{x\not= y} \fr{|f(x) - f(y)|}{|x-y|^\alpha}.
\end{equation}

We first derive a weighted gradient estimate for the heat kernel in Proposition \ref{heat_kernel_estimate}, inspired by \cite[Appendix 1]{constantin2016critical} and \cite[Appendix]{https://doi.org/10.1002/cpa.21720}.
\begin{prop}[Heat Kernel Estimate] Let $s 
\in (0,2)$ and $q \in [1, \infty).$ Then it holds that
\label{heat_kernel_estimate}
\begin{equation}
\label{heat_kernel_result}
    \int_0^\infty  t^{-1+\fr{s}{4}} \left( \int_\Omega |x-y|^{q\alpha} \, |(\na_x + \na_y) H_D (t,x,y)|^q \ dy \right)^\fr{1}{q} \ dt \leq C d(x)^{\alpha + \fr{s}{2} - d + \fr{d}{q} - 1},
\end{equation}    
for $-\frac{d}{q}< \alpha < 1 + d -\fr{d}{q} - \fr{s}{2}$.
\end{prop}
\noindent The proof of Proposition~\ref{heat_kernel_estimate} is presented in Appendix~\ref{kernel_appendix}.

Now we make use of Proposition \ref{heat_kernel_estimate} to prove the following pointwise commutator estimate:
\begin{thm}[Pointwise Commutator Estimate]
\label{commutator_theorem_1}   
Let $s \in (0,2)$. Let $q$ be a smooth scalar function and $v$ a smooth divergence-free vector field. Then it holds that 
\be\begin{aligned}
    \label{pointwise_estimate_result}
    \lvert \l^{-\fr{s}{2}} &(v \cdot \na q)(x) - (v \cdot \na \l^{-\fr{s}{2}} q)(x) \rvert \\
    & \lesssim [v]_{C^\alpha} \LnormD{\fr{q(y) - q(x)}{|x-y|^\beta}}{2}{dy} + [v]_{C^\gamma} \, |q(x)| \\
    &\hspace{2em} + |v(x)| \, \LnormD{\fr{q(y) - q(x)}{|x-y|^\delta}}{2}{dy} d(x)^{\delta + \fr{s}{2} - 2} + |v(x)| \, |q(x)| \, d(x)^{-1 + \fr{s}{2}},
\end{aligned}\ee
where $0 < \alpha, \gamma < 1$, $\alpha + \beta > 2 - \fr{s}{2}$, $\gamma > 1 - \fr{s}{2}$, and $0 \leq \delta < 2 - \fr{s}{2}$.
\end{thm}
\begin{proof}
Using the representation \eqref{negative_representation} for the negative fractional Laplacian, we study the pointwise difference of the following commutator and decompose it into the sum of two terms $\mathcal{A}$ and $\mathcal{B}$ 
\begin{equation}\begin{aligned}
\label{initial_representation}
    &\l^{-\fr{s}{2}} (v \cdot \na q)(x) - (v \cdot \na \l^{-\fr{s}{2}} q) (x) \\
    &= c_s \ComInt{H_D(x,y,t) \ \na_y \cdot (vq)} - c_s \ComInt{\na_x H_D(x,y,t) \cdot v(x) q(y)} \\
    &= - c_s \ComInt{\na_y H_D \cdot v(y) q(y)} - c_s \ComInt{\na_x H_D \cdot v(x) q(y)} \\
    &= -c_s \ComInt{\na_y H_D \cdot (v(y) - v(x)) q(y)} \\
    & \hspace{3em} - c_s \ComInt{(\na_x + \na_y) H_D \cdot v(x) q(y)} \\
    &\coloneqq \mathcal{A} + \mathcal{B},
\end{aligned}\end{equation}
by integrating the first term by parts and adding and subtracting the term $c_s \ComInt{\na_y H_D \cdot v(x) q(y)}$. Next, to handle the $\mathcal{A}$ term in \eqref{initial_representation}, we further add and subtract $c_s\ComInt{\na_y H_D \cdot (v(y) - v(x)) q(x)}$, which yields
\be\begin{aligned}
    \mathcal{A} &= - c_s \ComInt{\na_yH_D \cdot (v(y) - v(x)) \,(q(y)-q(x))} \\
    &\hspace{4em}- c_s \ComInt{\na_y H_D \cdot (v(y) - v(x)) \, q(x)} \\
    &\coloneqq \mathcal{A}_1 + \mathcal{A}_2.
\end{aligned}\ee
Now, we use the definition of the H\"older seminorm \eqref{holder_seminorm} and the gradient estimate \eqref{gradient_bound} to simplify the term $\mathcal{A}_1$ to
\be\begin{aligned}
    \left| \mathcal{A}_1 \right| &= c_s \left| \ComInt{\na_y H_D \cdot (v(y) - v(x)) \, (q(y) - q(x))} \right| \\
    &= c_s \left| \ComInt{|x-y|^{\alpha+\beta} \ \na_y H_D \cdot \fr{(v(y) - v(x))}{|x-y|^\alpha} \, \fr{(q(y) - q(x))}{|x-y|^\beta}} \right| \\
    &\lesssim [v]_{C^\alpha} \ComInt{|x-y|^{\alpha+\beta} \ |\na_y H_D| \ \left| \fr{q(y) - q(x)}{|x-y|^\beta} \right|}  \\
    &\lesssim [v]_{C^\alpha} \int_0^\infty t^{-\fr{5}{2}+\fr{s}{4}} \int_\Omega{e^{-\fr{|x-y|^2}{Kt}} \ |x-y|^{\alpha+\beta} \left| \fr{q(y) - q(x)}{|x-y|^\beta} \right|} \, dy \, dt.
\end{aligned}\ee
We then employ Fubini's theorem, the estimate \eqref{standard_integral_estimate}, and H\"older's inequality to reach 
\be\begin{aligned}
    \left| \mathcal{A}_1 \right| &\lesssim [v]_{C^\alpha} \int_\Omega \left( \int_0^\infty t^{-\fr{5}{2}+\fr{s}{4}}e^{-\fr{|x-y|^2}{Kt}} \, dt \right) \ |x-y|^{\alpha+\beta} \left| \fr{q(y) - q(x)}{|x-y|^\beta} \right| \, dy \\
    &\lesssim [v]_{C^\alpha} \int_\Omega |x-y|^{\alpha + \beta + \fr{s}{2}-3} \left| \fr{q(y) - q(x)}{|x-y|^\beta} \right| \, dy \\
    &\lesssim [v]_{C^\alpha} \LnormD{\fr{q(y) - q(x)}{|x-y|^\beta}}{2}{dy} \left( \int_\Omega |x-y|^{2\alpha + 2\beta + s -6}\, dy \right)^\fr{1}{2} \\
    &\lesssim [v]_{C^\alpha} \LnormD{\fr{q(y) - q(x)}{|x-y|^\beta}}{2}{dy},
\end{aligned}\ee
provided that $\alpha + \beta > 2 - \fr{s}{2}$ to ensure integrability.
$\mathcal{A}_2$ is estimated similarly through the definition of the H\"older seminorm \eqref{holder_seminorm}, the heat kernel bound \eqref{gradient_bound}, the standard estimate \eqref{standard_integral_estimate}, and the required condition for integrability, to obtain
\be\begin{aligned}
    \left| \mathcal{A}_2 \right| &= c_s \, | q(x)| \, \left| \ComInt{|x-y|^\gamma \ \na_y H_D \cdot \fr{(v(y) - v(x))}{|x-y|^\gamma}} \right| \\
    &\lesssim [v]_{C^\gamma} \, |q(x)| \ComInt{|x-y|^\gamma \ |\na_y H_D|} \\
    &\lesssim [v]_{C^\gamma} \, |q(x)| \int_\Omega \left( \int_0^\infty t^{-\fr{5}{2}+\fr{s}{4}} e^{-\fr{|x-y|^2}{Kt}} \right) |x-y|^\gamma \, dt \, dy \\
    &\lesssim [v]_{C^\gamma} \, |q(x)| \int_\Omega |x-y|^{\gamma + \fr{s}{2} - 3} \, dy \\
    &\lesssim [v]_{C^\gamma} \, |q(x)|,
\end{aligned}\ee
for a choice of $\gamma > 1 - \fr{s}{2}$. \\
Meanwhile, for the $\mathcal{B}$ term in \eqref{initial_representation}, we add and subtract $\ComInt{(\na_x + \na_y)H_D \cdot v(x)q(x)}$ to get the two treatable terms
\be\begin{aligned}
    \mathcal{B} &= -c_s \, v(x) \cdot \ComInt{[(\na_x + \na_y) H_D] (q(y) - q(x))} \\
    &\hspace{4em} - c_s \, v(x) \cdot \ComInt{[(\na_x + \na_y) H_D] q(x)} \\
    &\coloneqq \mathcal{B}_1 + \mathcal{B}_2.
\end{aligned}\ee
In view of H\"older's inequality and Proposition \ref{heat_kernel_estimate}, $\mathcal{B}_1$ simplifies to
\be\begin{aligned}
    | \mathcal{B}_1 | &= c_s \, |v(x)| \, \left| \ComInt{[|x-y|^\delta  \, (\na_x + \na_y) H_D] \fr{(q(y) - q(x))}{|x-y|^\delta}} \right| \\
    &\lesssim |v(x)| \, \LnormD{\fr{q(y) - q(x)}{|x-y|^\delta}}{2}{dy} \int_0^\infty t^{-1+\fr{s}{4}} \left( \int_\Omega |x-y|^{2\delta} \ |(\na_x + \na_y) H_D|^2 \, dy \right)^\fr{1}{2} \, dt \\
    &\lesssim |v(x)| \, \LnormD{\fr{q(y) - q(x)}{|x-y|^\delta}}{2}{dy} d(x)^{\delta + \fr{s}{2} - 2},
\end{aligned}\ee
where $\delta < 2 - \fr{s}{2}$.
Finally, $\mathcal{B}_2$ is estimated directly using Proposition \ref{heat_kernel_estimate}, yielding 
\be\begin{aligned}
    |\mathcal{B}_2| &= c_s \, |v(x)| \, |q(x)| \, \left| \ComInt{(\na_x + \na_y) H_D} \right| \\
    &\lesssim |v(x)| \, |q(x)| \, d(x)^{-1 + \fr{s}{2}}.
\end{aligned}\ee
Combining the estimates for all four terms, we reach the desired bound \eqref{pointwise_estimate_result}.
\end{proof}


With Theorem~\ref{commutator_theorem_1}, we further establish the following nonlinear commutator estimate.

\begin{thm}[Nonlinear Commutator Estimate] Let $\kappa > 0$. Let $s \in (0,2)$. Let $q \in L^2$ be a scalar function and $v \in H^2$ be a divergence-free vector field such that $v|_{\pa \Omega} = 0$.  Then it holds that 
\label{commutator_theorem_2}
\be\begin{aligned}
\label{nonlinear_estimate}
    &\OInt{\left| \left( \l^{-\fr{s}{2}} (v \cdot \na q) - v \cdot \na \l^{-\fr{s}{2}}q \right) \, \l^{-\fr{s}{2}}q \right|}  \leq \kappa \Lnorm{q}{2}^2 + C_{\kappa} \left(1 +  \Lnorm{\na v}{2}^\fr{2s}{2-s} \right) \Lnorm{\D v}{2}^2 \Lnorm{\l^{-\fr{s}{2}}q}{2}^2.
\end{aligned}\ee
\end{thm}
\begin{proof}  The following proof involves boundedness of the linear-in-$q$ operator $[\Lambda^{-\frac{s}{2}},v \cdot \na]q$ from $L^2$ into some $L^p$ spaces for a fixed $v \in H^2$ with vanishing boundary condition. The latter is established via application of Theorem \ref{commutator_theorem_1} for $v \in H^2$ and $q \in C^{\infty}(\Omega)$ followed by a classical extension argument that holds due to the density of $C^{\infty}(\Omega)$ in $L^2$. This technicality will be disregarded for simplicity.

Fix $\varepsilon > 0$ to be determined later. We study the commutator $[\l^{-\fr{s}{2}}, v \cdot \na]q$ multiplied by $\l^{-\fr{s}{2}}q$ in $L^2$ by employing the pointwise bound derived in Theorem \ref{commutator_theorem_1} (with $\delta = 1 - \varepsilon$) to get 
\be\begin{aligned}
    &\OInt{\left| \left(\l^{-\fr{s}{2}} (v \cdot \na q) - v \cdot \na \l^{-\fr{s}{2}} q \right) \, \l^{-\fr{s}{2}}q \right|} \\
    &\hspace{2em}\lesssim \OInt{\left([v]_{C^\alpha} \LnormD{\fr{q(y)-q(x)}{|x-y|^\beta}}{2}{dy} |\l^{-\fr{s}{2}} q(x)| + [v]_{C^\gamma} |q(x)| |\l^{-\fr{s}{2}} q(x)|\right)} \\
    &\hspace{4em}+ \OInt{\fr{|v(x)|}{d(x)^{1 - \fr{s}{2} + \varepsilon}} \LnormD{\fr{q(y) - q(x)}{|x-y|^{\delta}}}{2}{dy} |\l^{-\fr{s}{2}} q (x)|} + \OInt{\fr{|v(x)|}{d(x)^{1-\fr{s}{2}}} \, |q(x)| \, |\l^{-\fr{s}{2}}q (x)|}.
\end{aligned}\ee
We define
\begin{align}
    \mathcal{M} &\coloneqq \OInt{ \left( [v]_{C^\alpha} \LnormD{\fr{q(y)-q(x)}{|x-y|^\beta}}{2}{dy} |\l^{-\fr{s}{2}} q(x)| + [v]_{C^\gamma} |q(x)| \,|\l^{-\fr{s}{2}} q(x)| \right)}, \\
    \mathcal{N} &\coloneqq \OInt{\fr{|v(x)|}{d(x)^{1 - \fr{s}{2} + \varepsilon}} \LnormD{\fr{q(y) - q(x)}{|x-y|^{\delta}}}{2}{dy} |\l^{-\fr{s}{2}} q (x)|},
\end{align}
and
\begin{equation}
    \mathcal{O} \coloneqq \OInt{\fr{|v(x)|}{d(x)^{1-\fr{s}{2}}} \, |q(x)| \, |\l^{-\fr{s}{2}}q (x)|}.
\end{equation}
For $\mathcal{M}$, H\"older's inequality yields
\be\begin{aligned}
    \mathcal{M} \leq [v]_{C^\alpha} \Lnorm{\l^{-\fr{s}{2}}q}{2} \left( \int_\Omega \int_\Omega \fr{|q(y) - q(x)|^2}{|x-y|^{2\beta}} \, dy \, dx \right)^\fr{1}{2} + [v]_{C^\gamma} \Lnorm{q}{2} \Lnorm{\l^{-\fr{s}{2}}q}{2},
\end{aligned}\ee
where we then use the identity $(q(y)-q(x))^2 \lesssim q(y)^2 + q(x)^2$ and Fubini's theorem to get that
\be\begin{aligned}
\label{q_integral_bound}
    \int_\Omega \int_\Omega \fr{|q(y) - q(x)|^2}{|x-y|^{2\beta}} \, dy \, dx &\lesssim \int_\Omega \int_\Omega \fr{q(y)^2 + q(x)^2}{|x-y|^{2\beta}} \, dy \, dx \\
    &= C \int_\Omega q(y)^2 \int_\Omega \fr{1}{|x-y|^{2\beta}} \, dx \, dy + C \int_\Omega q(x)^2 \int_\Omega \fr{1}{|x-y|^{2\beta}} \, dy \, dx \\
    &\lesssim \LnormD{q}{2}{\Omega}^2,
\end{aligned}\ee
by choosing $\beta < 1$ to ensure the integrability of the inner integrals. Furthermore, for $m < 1$, as $v$ vanishes on the boundary, the following embedding holds
\begin{equation}
    [v]_{C^m} \lesssim \Lnorm{\D v}{2}.
\end{equation}
Hence, in view of the above bounds and Young's inequality, the term $\mathcal{M}$ is estimated as
\begin{equation}
    \mathcal{M} \lesssim \Lnorm{q}{2} \Lnorm{\D v}{2} \Lnorm{\l^{-\fr{s}{2}} q}{2} \leq \fr{\kappa}{3} \Lnorm{q}{2}^2 + C \Lnorm{\D v}{2}^2 \Lnorm{\l^{-\fr{s}{2}}q}{2}^2.
\end{equation}
Moving to the term $\mathcal{N}$, we apply the estimate \eqref{q_integral_bound} (provided that $\delta < 1$) and the H\"older, Gagliardo–Nirenberg, and Young inequalities, obtaining
\be\begin{aligned}
    \mathcal{N} &\lesssim \Lnorm{q}{2} \, \Lnorm{\l^{-\fr{s}{2}}q}{\fr{8}{4-s^2}} \, \Lnorm{\fr{v}{d^{1-\fr{s}{2}+\varepsilon}}}{\fr{8}{s^2}} \lesssim \Lnorm{q}{2}^{1+\fr{s}{2}} \, \Lnorm{\l^{-\fr{s}{2}}q}{2}^{1-\fr{s}{2}} \, \Lnorm{\fr{v}{d^{1-\fr{s}{2} + \varepsilon}}}{\fr{8}{s^2}} \\
    &\leq \fr{\kappa}{3} \Lnorm{q}{2}^2 + C \Lnorm{\l^{-\fr{s}{2}}q}{2}^2 \, \Lnorm{\fr{v}{d^{1-\fr{s}{2}+\varepsilon}}}{\fr{8}{s^2}}^\fr{4}{2-s}.
\end{aligned}\ee
We now deal with the last factor by appealing to Rellich's inequality (see \cite{edmunds2005weighted} and references therein),
the chain rule (where $|\na(|v|^\alpha)| \leq \alpha |v|^{\alpha-1} |\na v|$), and H\"older's inequality in the following fashion  
\be\begin{aligned}
    \Lnorm{\fr{v}{d^{1-\fr{s}{2}+\varepsilon}}}{\fr{8}{s^2}}^\fr{4}{2-s} &= \left( \OInt{\left(\fr{|v(x)|^{\fr{1}{1-\fr{s}{2} + \varepsilon}}}{d(x)}\right)^{\fr{8-4s+8\varepsilon}{s^2}}} \right)^\fr{s^2}{2(2-s)} \lesssim \Lnorm{\na \left(|v|^\fr{1}{1-\fr{s}{2} + \varepsilon} \right)}{\fr{8-4s+8\varepsilon}{s^2}}^{\fr{4-2s+4\varepsilon}{2-s}} \\
    &\lesssim \Lnorm{\na v \ |v|^{\fr{s - 2\varepsilon}{2-s+2\varepsilon}}}{\fr{8 - 4s+8\varepsilon}{s^2}}^\fr{4-2s+4\varepsilon}{2-s} \lesssim \Lnorm{\na v}{\fr{2-s+2\varepsilon}{\varepsilon}}^\fr{4-2s+4\varepsilon}{2-s} \Lnorm{|v|^\fr{s-2\varepsilon}{2-s+2\varepsilon}}{\fr{8-4s+8\varepsilon}{s^2-4\varepsilon}}^\fr{4-2s+4\varepsilon}{2-s} \\
    &= C \Lnorm{\na v}{\fr{2-s+2\varepsilon}{\varepsilon}}^\fr{4-2s+4\varepsilon}{2-s} \Lnorm{v}{\fr{4(s-2\varepsilon)}{s^2-4\varepsilon}}^\fr{2(s-2\varepsilon)}{2-s} \\
\end{aligned}\ee
provided that $\varepsilon \leq \fr{s}{2}$, $\varepsilon < \fr{s^2}{4}$, and $\varepsilon \geq \fr{s^2}{8} + \fr{s}{2} - 1$. These three conditions are satisfied simultaneously whenever $0 < s < 2$. We note that the first condition is required for the proper differentiation of
$\nabla \left(|v|^{\frac{1}{1-\frac{s}{2}+\epsilon}}\right)$,
ensuring that only nonnegative powers of $|v|$ appear. The second condition is
necessary for the application of H\"older's inequality. Finally, the last
condition allows us to use Rellich's inequality, which holds in $L^p$
spaces for $p \ge 1$.
Next, we use the Gagliardo–Nirenberg interpolation inequality to get
\be\begin{aligned}
    \Lnorm{\fr{v}{d^{1-\fr{s}{2}+\varepsilon}}}{\fr{8}{s^2}}^\fr{4}{2-s} &\lesssim \Lnorm{\na v}{2}^\fr{4\varepsilon}{2-s} \Lnorm{\D v}{2}^2 \Lnorm{v}{\fr{4(s-2\varepsilon)}{s^2-4\varepsilon}}^\fr{2(s-2\varepsilon)}{2-s}.
\end{aligned}\ee
Using classical continuous Sobolev embeddings, the latter boils down to
\begin{equation}
    \Lnorm{\fr{v}{d^{1-\fr{s}{2}+\varepsilon}}}{\fr{8}{s^2}}^\fr{4}{2-s} \lesssim \Lnorm{\na v}{2}^\fr{2s}{2-s} \Lnorm{\D v}{2}^2.
\end{equation}
Finally, the term $\mathcal{O}$ is similarly estimated as 
\be\begin{aligned}
    \mathcal{O} &\lesssim \Lnorm{q}{2} \, \Lnorm{\l^{-\fr{s}{2}}q}{\fr{8}{4-s^2}} \, \Lnorm{\fr{v}{d^{1 - \fr{s}{2}}}}{\fr{8}{s^2}} \lesssim \Lnorm{q}{2}^{1 + \fr{s}{2}} \, \Lnorm{\l^{-\fr{s}{2}}q}{2}^{1 - \fr{s}{2}} \, \Lnorm{\fr{v}{d^{1-\fr{s}{2}}}}{\fr{8}{s^2}} \\
    &\leq \fr{\kappa}{3} \Lnorm{q}{2}^2 + C \Lnorm{\l^{-\fr{s}{2}}q}{2}^2 \, \Lnorm{\fr{v}{d^{1-\fr{s}{2}}}}{\fr{8}{s^2}}^{\fr{4}{2-s}} \leq \fr{\kappa}{3} \Lnorm{q}{2}^2 + C \Lnorm{\l^{-\fr{s}{2}}q}{2}^2 \, \Lnorm{\fr{v}{d^{1-\fr{s}{2}+\varepsilon}}}{\fr{8}{s^2}}^{\fr{4}{2-s}},
\end{aligned}\ee
where the last inequality follows from the boundedness of $\Omega$.
We proceed as for the term $\mathcal{N}$ and combine these estimates, establishing the desired bound \eqref{nonlinear_estimate}.
\end{proof}

\section{Existence and Uniqueness of Solutions}
\label{sec_existence_and_uniqueness}
Our main goal is to prove the global well-posedness of the magnetically forced electroconvection system \eqref{main-system} for $q_0 \in L^4$ and $u_0 \in V$. 
To this end, we consider the following parabolically $\eta$-regularized system 
\begin{subequations}
\begin{align}
    &\partial_t u^\eta + \PP (u^\eta \cdot \na u^\eta) + \nu A u^\eta + B^2 u^\eta + B \, \PP(q^\eta (u^\eta)^\perp) = B R^\perp q^\eta - \PP  (q^\eta R q^\eta), \label{m_navier_stokes_eq} \\
    &\partial_t q^\eta + u^\eta \cdot \na q^\eta + \l q^\eta - B \Div{(u^\eta)^\perp} - \eta \D q^\eta = 0 \label{m_charge_density_eq}, \\
    &\Div{u^\eta} = 0,
\end{align}
\end{subequations}
{with initial conditions $u_0^\eta = u_0$ and $q_0^\eta = J_\eta q_0$, and boundary conditions $u^\eta |_{\partial \Omega} = q^\eta|_{\partial \Omega} = 0$}. Here $\PP$ denotes the Leray projector, $A = -\PP \Delta$ is the Stokes operator, and $J_\eta$ is the spectrally regularizing mollifier, defined as 
\begin{equation}
\label{mollifier_definition}
    J_\eta f(x) \coloneqq - \fr{1}{\ln \eta} \int_\eta^\fr{1}{\eta} \fr{e^{t \D_D} f(x)}{t} \ dt,
\end{equation}
for $0< \eta < 1$. This family of regularizers smooths out any $L^2$ function in the following sense: given a function $\theta \in L^2(\Omega)$, we have $J_{\eta}\theta \in \mathcal{D}(\l^k)$ for any $k \in \mathbb{N}$. In addition, $J_{\eta}$ is uniformly bounded on $\mathcal{D}(\l^k)$ for any $k \in \mathbb{Z}$ and on $L^p$ spaces for any $p \in [1,\infty]$. Also, for $f \in \mathcal{D}(\l^k)$, we have $J_\eta f \to f$ in $\mathcal{D}(\l^k)$ as $\eta \to 0$ (see \cite{abdo2024regularity, abdo2024dirichlet} for more details). 

Without loss of generality, we consider the magnetic constant $B$ to be nonnegative.



\begin{thm}
   Let $q_0 \in L^4$ and $u_0 \in V$. Let $T > 0$ be arbitrary. Then system \eqref{main-system} has a unique solution $(q,u)$ on $[0, T]$ with regularity 
\begin{equation} 
q \in L^{\infty}(0,T; L^4) \cap L^2(0,T; \mathcal{D} (\l^{\frac{1}{2}})), \quad u \in L^{\infty}(0,T; V) \cap L^2(0,T; \mathcal{D}(A)).
\end{equation} Moreover, there is a constant $C_B$ depending only on $B$, the $L^4$ norm of the initial data, $H^1$ norm of the initial velocity, and the parameters of the problem, and a positive universal constant $c_0$ independent of $B$, such that 
\be \label{decay-thm}
\|u(t)\|_{H^1} + \|q(t)\|_{L^4} \le C_Be^{- c_0\min \left( 1, \fr{1}{B^2 + \nu} \right) t}
\ee for any $t \ge 0$. 
\end{thm}
\begin{proof}
For each $\eta>0$, the mollified system is parabolically regularized and therefore admits a unique solution $(q^\eta,u^\eta)$, with
\[
q^\eta \in L^\infty(0,T;H^1)\cap L^2(0,T;H^2),
\]
and
\[
u^\eta \in L^\infty(0,T;\mathcal D(A^{1/2}))\cap L^2(0,T;\mathcal D(A)).
\]
The latter can be obtained through a suitable Galerkin approximation scheme (see, for instance, \cite{abdo2025electroconvection, abdo2023long}). However, Galerkin approximations are not well suited for deriving $L^4$-estimates, since they do not preserve the transport structure associated with divergence-free vector fields. Consequently, $H^1$-estimates for $q$ must be established, which justifies mollifying the initial charge density, assumed only to belong to $L^4$.

Having constructed these smooth solutions, we derive estimates uniform in $\eta$ and then pass to the limit by means of the Aubin--Lions and Banach--Alaoglu theorems, and the lower semicontinuity of norms. This is based on a standard compactness argument, and the details will be omitted. We present below the nontrivial uniform-in-$\eta$ decay estimates, which rely crucially on the coupled energy structure of the system.

\medskip

{\bf{Step 1.} Decay for $(q^{\eta},u^{\eta})$ in $\mathcal{D}(\l^{-\fr{1}{2}}) \times L^2$.} Multiplying \eqref{m_navier_stokes_eq} by $u^\eta$ and integrating spatially, we have
\begin{equation}
    \fr{1}{2} \dt \Lnorm{u^\eta}{2}^2 + \nu \Lnorm{A^\fr{1}{2} u^\eta}{2}^2 + B^2 \Lnorm{u^\eta}{2}^2 = B \OInt{R^\perp q^\eta \cdot u^\eta} - \OInt{ (q^\eta R q^\eta) \cdot u^\eta}, \label{i_navier_stokes_eq}
\end{equation}
wherein we used the divergence-free condition on $u^\eta$, the self-adjointness of $\PP$, and the fact that $(u^\eta)^\perp \cdot u^\eta = 0$ to obtain 
\[
\OInt{\PP(u^\eta \cdot \na u^\eta) \cdot u^\eta}=0, \quad \OInt{\PP (q^\eta (u^\eta)^\perp) \cdot u^\eta} = 0.
\]
Meanwhile, multiplying the charge density equation \eqref{m_charge_density_eq} by $\l^{-1} q^\eta$, we obtain the equation for the $L^2$ evolution of $\l^{-\fr{1}{2}} q^\eta$, 
\begin{equation}
    \fr{1}{2} \dt \Lnorm{\l^{-\fr{1}{2}} q^\eta}{2}^2 + \eta \Lnorm{\na \l^{-\fr{1}{2}} q^\eta}{2}^2 + \Lnorm{q^\eta}{2}^2 = - \OInt{u^\eta \cdot \na q^\eta \l^{-1} q^\eta} + B \OInt{\Div{(u^\eta)^\perp} \ \l^{-1} q^\eta}. \label{i_charge_density_eq}
\end{equation}
Using the divergence-free property of $u^\eta$,  integration by parts, and the homogeneous Dirichlet boundary condition obeyed by $u^{\eta}$, we have 
\be\begin{aligned}
    - \OInt{u^\eta \cdot \na q^\eta \ \l^{-1} q^\eta} &= - \OInt{\Div{(u^\eta q^\eta)} \ \l^{-1} q^\eta} = \OInt{u^\eta q^\eta \cdot \na \l^{-1} q^\eta} = \OInt{u^\eta \cdot ( q^\eta R q^\eta) }  .
\end{aligned}\ee
Now, adding \eqref{i_navier_stokes_eq} and \eqref{i_charge_density_eq}, cancelling the nonlinear advective terms, and dropping the nonnegative regularization term from the left-hand side, we obtain the coupled differential inequality
\be\begin{aligned}
    \fr{1}{2} &\dt \left(\Lnorm{u^\eta}{2}^2 + \Lnorm{\l^{-\fr{1}{2}}q^\eta}{2}^2 \right) + \nu \Lnorm{A^\fr{1}{2} u^\eta}{2}^2 + B^2 \Lnorm{u^\eta}{2}^2 + \Lnorm{q^\eta}{2}^2 \leq 2B \OInt{R^\perp q^\eta \cdot u^\eta},
    \label{diff_inequality_first_step}
\end{aligned}\ee
where we integrate by parts and use the fact that rotation by $\frac{\pi}{2}$ is skew-symmetric to get
\begin{equation}
    \OInt{\na \cdot(u^\eta)^\perp \ \l^{-1} q^\eta} = - \OInt{(u^\eta)^\perp \cdot Rq^\eta} = \OInt{u^\eta \cdot R^\perp q^\eta}.
\end{equation}
Upon application of the H\"older and Young inequalities, and the boundedness of the Riesz transform in $L^2$ (via the identity $\Lnorm{\na \l^{-1}q}{2}^2 = \Lnorm{\l \l^{-1} q}{2}^2= \Lnorm{q}{2}^2$), the right-hand side of \eqref{diff_inequality_first_step} becomes
\begin{equation}
\label{est1}
    2B \OInt{R^\perp q^\eta \cdot u^\eta} \le 2B \Lnorm{q^\eta}{2} \Lnorm{u^\eta}{2} \leq \fr{B^2}{B^2 + \delta} \Lnorm{q^\eta}{2}^2 + (B^2+\delta) \Lnorm{u^\eta}{2}^2.
\end{equation}
Moving these terms to the dissipation of \eqref{diff_inequality_first_step}, splitting the $\nu \Lnorm{A^\fr12 u^\eta}{2}^2$ into two, and applying the Poincar\'e inequality on $C \Lnorm{u^\eta}{2}^2 \leq \Lnorm{\na u^\eta}{2}^2 = \Lnorm{A^\fr{1}{2} u^\eta}{2}^2$ and $C \Lnorm{\l^{-\fr{1}{2}} q^\eta}{2}^2 \leq \Lnorm{q^\eta}{2}^2$, we reach
\begin{equation}
    \fr12 \dt \left(\Lnorm{u^\eta}{2}^2 + \Lnorm{\l^{-\fr{1}{2}}q^\eta}{2}^2 \right) + \left(\fr{C {\nu}}{2} - \delta\right) \Lnorm{u^\eta}{2}^2 + \fr{\nu}{2} \Lnorm{A^\fr{1}{2} u^\eta}{2}^2+ \frac{C \delta}{B^2+\delta}\Lnorm{\l^{- \frac{1}{2}}q^\eta}{2}^2 \leq 0.
\end{equation}
Now, we can choose $\delta = \fr{C\nu}{4}$ to get the inequality 
\begin{equation}
    \dt \left( \Lnorm{u^\eta}{2}^2 + \Lnorm{\l^{-\fr{1}{2}} q^\eta}{2}^2 \right) + \fr{C \nu}{2} \Lnorm{u^\eta}{2}^2 + \fr{2C^2 \nu}{4B^2 + C\nu} \Lnorm{\l^{-\fr{1}{2}} q^\eta}{2}^2 + \nu \Lnorm{A^\fr{1}{2} u^\eta}{2}^2 \leq 0,
\end{equation}
which can be rewritten as 
\begin{equation}
    \dt \left( \Lnorm{u^\eta}{2}^2 + \Lnorm{\l^{-\fr{1}{2}} q^\eta}{2}^2 \right) + {C_\nu} \min \left( 1, \fr{1}{B^2 + \nu} \right) \left( \Lnorm{u^\eta}{2}^2 + \Lnorm{\l^{-\fr{1}{2}} q^\eta}{2}^2 \right) + \nu \Lnorm{A^\fr{1}{2} u^\eta}{2}^2 \leq 0.
\end{equation} 
We now apply the Gr\"onwall Lemma~\ref{gronwall_lemma} on the latter differential inequality with $c_1 = {C_\nu}\min \left( 1, \fr{1}{B^2 + \nu}\right), \, c_2 = \nu, \, F = G = H =I= 0$, and $a = b = 2c_1$, obtaining 
\begin{equation}
\label{decay_of_u_and_lambda_minus_half_q_and_time_integrability_of_A_half_u}
    \Lnorm{u^\eta}{2}^2 + \Lnorm{\l^{-\fr{1}{2}} q^\eta}{2}^2 \leq \G_1 e^{- C_1 \min \left( 1, \fr{1}{B^2 + \nu}\right) t}, \quad \text{ and } \quad \TInt{\Lnorm{A^\fr{1}{2} u^\eta}{2}^2} \leq \G_2,
\end{equation}
where $\G_1$ and $\G_2$ are constants depending on $\Lnorm{u_0}{2}^2$, $\Lnorm{\l^{-\fr{1}{2}} q_0}{2}^2$ but are independent of $B$.

\medskip

{\bf{Step 2.} Decay of the regularized charge densities in $L^2$.} Multiplying \eqref{m_charge_density_eq} by $q^\eta$, integrating spatially over $\Omega$, and using the divergence-free property of $u^\eta$, we have
\be\begin{aligned}
\label{step2_diff_inequality_start}
    \fr{1}{2} \dt \Lnorm{q^\eta}{2}^2 + \Lnorm{\l^\fr{1}{2} q^\eta}{2}^2 + \eta \Lnorm{\na q^\eta}{2}^2 = B\OInt{\Div{(u^\eta)^\perp} \ q^\eta}.
\end{aligned}\ee
We estimate the term on the right-hand side by introducing $\l^{-\fr{1}{2}}$, $\l^\fr{1}{2}$, interpolating $\l^\fr{1}{2} u^\eta$ between $u^\eta$ and $A^\fr{1}{2} u^\eta$ in $L^2$, and applying H\"older's and Young's inequalities to get

\be\begin{aligned}
\label{step2_estimate_1_term}
    B\OInt{\Div{ (u^\eta)^\perp \ q^\eta}} &= B\OInt{\l^{-\fr{1}{2}} \Div{ (u^\eta)^\perp \ \l^\fr{1}{2} q^\eta}} \leq C B \Lnorm{\l^\fr{1}{2} u^\eta}{2} \Lnorm{\l^\fr{1}{2} q^\eta}{2} \\
    &\leq \fr{1}{2} \Lnorm{\l^\fr{1}{2} q^\eta}{2}^2 + CB^2 \Lnorm{\l^\fr{1}{2} u^\eta}{2}^2 \\
    &\leq \fr{1}{2} \Lnorm{\l^\fr{1}{2} q^\eta}{2}^2 + CB^2  \Lnorm{u^\eta}{2} \Lnorm{A^\fr{1}{2} u^\eta}{2} \\
    &\leq \fr{1}{2} \Lnorm{\l^\fr{1}{2} q^\eta}{2}^2 + CB^2  \Lnorm{u^\eta}{2} \left(1+\Lnorm{A^\fr{1}{2} u^\eta}{2}^2\right).
\end{aligned}\ee
Combining \eqref{step2_estimate_1_term} with the inequality in \eqref{step2_diff_inequality_start} and applying the Poincar\'e inequality on half of the dissipation, we obtain 
\be\begin{aligned}
   \dt \Lnorm{q^\eta}{2}^2 + C_2 \Lnorm{ q^\eta}{2}^2 + \fr{1}{2} \Lnorm{\l^\fr{1}{2} q^\eta}{2}^2 \leq CB^2 \left( \Lnorm{A^\fr{1}{2} u^\eta}{2}^2 + 1 \right) \, \Lnorm{u^\eta}{2}. 
\end{aligned}\ee
We again make use of the Gr\"onwall Lemma~\ref{gronwall_lemma} with $c_1 = C_2, \, c_2 = 1/2, \, F = 0$,
\begin{equation}
    G(t) = CB^2 \Lnorm{A^\fr{1}{2} u^\eta}{2}^2, \quad H(t) =\Lnorm{u^\eta}{2}, \quad I(t) = CB^2 \Lnorm{u^\eta}{2}, \quad a = b = C_1 \min \left( 1, \fr{1}{B^2 + \nu} \right),
\end{equation}
where the time integrability of $G$ and the decay rates of $H$ and $I$ are justified by \eqref{decay_of_u_and_lambda_minus_half_q_and_time_integrability_of_A_half_u}. Hence, we reach
\begin{equation}
\label{decay_of_q_and_lambda_half_q_and_time_integrability}
    \Lnorm{q^\eta}{2}^2 \leq M \G_3 e^{- C_3\min \left(1, \fr{1}{B^2 + \nu} \right) t},
    \quad \text{ and } \quad 
    \TInt{\Lnorm{\l^\fr{1}{2} q^\eta}{2}^2} \leq M \G_4,
\end{equation}
where $M = \left( 1 + B^2 + \fr{B^2}{\min \left(1 , \fr{1}{B^2 + \nu} \right)} \right)$, and $\G_3, \, \G_4$ are constants depending on the parameters (except $B$) and on $\Lnorm{q_0}{2}$ and $\|u_0\|_{L^2}$.

\medskip

{\textbf{Step 3. Decay of $(q^{\eta}, u^{\eta})$ in $L^4 \times V$.}} We take the scalar product in $L^2$ of the equation \eqref{m_navier_stokes_eq} with $A u^\eta$ and integrate spatially, reaching
\be\begin{aligned}
\label{H1_u_eps_first_step}
   \fr{1}{2} \dt &\Lnorm{A^\fr{1}{2} u^\eta}{2}^2 +  B^2 \Lnorm{A^\fr{1}{2} u^\eta}{2}^2 + \nu \Lnorm{A u^\eta}{2}^2 \\
   &= - \OInt{\PP (u^\eta \cdot \na u^\eta) Au^\eta}
   - \OInt{\PP  (q^\eta R q^\eta) A u^\eta} - B \OInt{\PP(q^\eta (u^\eta)^\perp) A u^\eta} \\
   & \quad \quad + B \OInt{R^\perp q^\eta Au^\eta} \\
   &= \text{I} + \text{II} + \text{III} + \text{IV}.
\end{aligned}\ee
We proceed by treating the terms individually. For the convective nonlinear term I, we use the self-adjointness and ellipticity of the Stokes operator, the Ladyzhenskaya interpolation inequality, and Young's inequality, to obtain
\be\begin{aligned}
    |\text{I}| &= \left| \OInt{\PP (u^\eta \cdot \na u^\eta) A u^\eta} \right| \leq \Lnorm{u^\eta \cdot \na u^\eta}{2} \Lnorm{A u^\eta}{2} \leq C \Lnorm{u^\eta}{4} \Lnorm{\na u^\eta}{4} \Lnorm{A u^\eta}{2} \\
    &\leq C \Lnorm{u^\eta}{2}^\fr{1}{2} \Lnorm{\na u^\eta}{2} \Lnorm{A u^\eta}{2}^\fr{3}{2} \leq \fr{\nu}{8} \Lnorm{A u^\eta}{2}^2 + C\Lnorm{u^\eta}{2}^2 \Lnorm{\na u^\eta}{2}^4.
\end{aligned}\ee
By the boundedness of the Riesz transform in $L^4$, Young's inequality, and the interpolation inequality of $L^p$ norms, we estimate the second term II as follows
\be\begin{aligned}
\label{oo}
    |\text{II}| &= \left| \OInt{\PP  (q^\eta R q^\eta) A u^\eta} \right| \leq C \Lnorm{q^\eta Rq^\eta}{2} \Lnorm{A u^\eta}{2} \leq C \Lnorm{q^\eta}{4}^2 \Lnorm{A u^\eta}{2} \\
    &\leq \fr{\nu}{8} \ \Lnorm{A u^\eta}{2}^2 + C \Lnorm{q^\eta}{4}^4 \leq \fr{\nu}{8} \Lnorm{A u^\eta}{2}^2 + C \Lnorm{q^\eta}{8}^\fr{8}{3} \Lnorm{q^\eta}{2}^\fr{4}{3}.
\end{aligned}\ee
Using the embedding of $\mathcal{D}(\l^\fr{1}{2})$ in $L^4$, we have
\begin{equation}
\label{dissipation_bound}
    \Lnorm{q^\eta}{8}^2 = \Lnorm{(q^\eta)^2}{4} \leq C \Lnorm{\l^\fr{1}{2} (q^\eta)^2}{2}.
\end{equation} 
Consequently, \eqref{oo} is bounded by
\begin{equation}
    |\text{II}| \leq \fr{\nu}{8} \Lnorm{A u^\eta}{2}^2 + C \Lnorm{\l^\fr{1}{2} (q^\eta)^2}{2}^\fr{4}{3} \Lnorm{q^\eta}{2}^\fr{4}{3}.
\end{equation}
Similarly, the third term is estimated, using the embedding of $H^\fr{1}{2}$ in $L^4$ and $H^1$ in $L^4$, as
\be\begin{aligned}
    |\text{III}| &= B \left| \OInt{\PP(q^\eta (u^\eta)^\perp) A u^\eta} \right| \leq B \Lnorm{q^\eta}{4} \Lnorm{u^\eta}{4} \Lnorm{A u^\eta}{2} \leq C B \Lnorm{\l^\fr{1}{2} q^\eta}{2} \Lnorm{\na u^\eta}{2} \Lnorm{A u^\eta}{2} \\
    &\leq \fr{\nu}{8} \Lnorm{A u^\eta}{2}^2 + C B^2 \Lnorm{\l^\fr{1}{2} q^\eta}{2}^2 \Lnorm{\na u^\eta}{2}^2.
\end{aligned}\ee
As for the fourth term IV, we get
\be\begin{aligned}
    |\text{IV}| = B \left| \OInt{R^\perp q^\eta Au^\eta} \right| \leq \fr{\nu}{8} \Lnorm{A u^\eta}{2}^2 + C B^2 \Lnorm{q^\eta}{2}^2.
\end{aligned}\ee
Combining the estimates of these four terms into \eqref{H1_u_eps_first_step}, we get
\be\begin{aligned}
\label{H1_velocity_final}
    \dt &\Lnorm{A^\fr{1}{2} u^\eta}{2}^2 + B^2 \Lnorm{A^\fr{1}{2} u^\eta}{2}^2 + \nu \Lnorm{A u^\eta}{2}^2 \\
    &\quad \leq C \Lnorm{u^\eta}{2}^2 \Lnorm{\na u^\eta}{2}^4 + C \Lnorm{\l^\fr{1}{2} (q^\eta)^2}{2}^\fr{4}{3} \Lnorm{q^\eta}{2}^\fr{4}{3} + C B^2 \Lnorm{\l^\fr{1}{2} q^\eta}{2}^2 \Lnorm{\na u^\eta}{2}^2 + C B^2 \Lnorm{q^\eta}{2}^2.
\end{aligned}\ee

We couple the $\mathcal{D} (A^\fr{1}{2})$ evolution of $u^\eta$ with the $L^4$ evolution of $q^\eta$. The latter can be obtained by multiplying equation \eqref{m_charge_density_eq} by $(q^\eta)^3$, which yields
\be\begin{aligned}
\label{L4_charge_first_inequality}
    \fr{1}{4} \dt \Lnorm{q^\eta}{4}^4 + C_4 \Lnorm{q^\eta}{4}^4 + C_5 \Lnorm{\l^\fr{1}{2} (q^\eta)^2}{2}^2
    \leq B \OInt{ \Div{(u^\eta)^\perp} (q^\eta)^3}, 
\end{aligned}\ee
where we have used the C\'ordoba-C\'ordoba inequality (see \cite{abdo2024dirichlet, cordoba2004maximum}) and the fact that the nonlinear term $\OInt{u^\eta \cdot \na q^\eta (q^\eta)^3}$ vanishes. Using the H\"older and Young inequalities, and the embedding of $\mathcal{D}(\l^\fr{1}{2})$ in $L^4$, we estimate the right-hand side as
\be\begin{aligned}
    B \OInt{ \Div{(u^\eta)^\perp (q^\eta)^3}} &\leq B \Lnorm{\na u^\eta}{2} \Lnorm{q^\eta}{6}^3 \leq B \Lnorm{\na u^\eta}{2} \Lnorm{q^\eta}{4} \Lnorm{q^\eta}{8}^2 \\
    &\leq B \Lnorm{\na u^\eta}{2} \Lnorm{\l^\fr{1}{2} q^\eta}{2} \Lnorm{q^\eta}{8}^2 \leq C B \Lnorm{\na u^\eta}{2} \Lnorm{\l^\fr{1}{2} q^\eta}{2} \Lnorm{\l^\fr{1}{2} (q^\eta)^2}{2} \\
    &\leq \fr{C_5}{2} \Lnorm{\l^\fr{1}{2} (q^\eta)^2}{2}^2 + C B^2 \Lnorm{\na u^\eta}{2}^2 
    \Lnorm{\l^\fr{1}{2} q^\eta}{2}^2,
\end{aligned}\ee
where we used \eqref{dissipation_bound} again.
Therefore, we obtain
\begin{equation}
\label{oiga}
    \dt \Lnorm{q^\eta}{4}^4 + 4 C_4 \Lnorm{q^\eta}{4}^4 + 2C_5 \Lnorm{\l^\fr{1}{2} (q^\eta)^2}{2}^2 \leq C B^2 \Lnorm{\na u^\eta}{2}^2 \Lnorm{\l^\fr{1}{2} q^\eta}{2}^2.
\end{equation}
Adding together \eqref{H1_velocity_final} and \eqref{oiga} and using Young's inequality again, we obtain
\be\begin{aligned}
    &\dt \left( \Lnorm{A^\fr{1}{2} u^\eta}{2}^2 + \Lnorm{q^\eta}{4}^4 \right) + \left( C_6 \nu \Lnorm{A^\fr{1}{2} u^\eta}{2}^2 + 4 C_4 \Lnorm{q^\eta}{4}^4 \right) + \left( \fr{\nu}{2} \Lnorm{A u^\eta}{2}^2 + C_5 \Lnorm{\l^\fr{1}{2} (q^\eta)^2}{2}^2 \right) 
    \\&\leq C \left( \Lnorm{u^\eta}{2}^2 \Lnorm{A^\fr{1}{2} u^\eta}{2}^2 + (1 + B^2) \Lnorm{\l^\fr{1}{2} q^\eta}{2}^2 \right) \left(\Lnorm{A^\fr{1}{2} u^\eta}{2}^2 +\|q^{\eta}\|_{L^4}^4\right) + C \Lnorm{q^\eta}{2}^4 + C B^2 \Lnorm{q^\eta}{2}^2,
\end{aligned}\ee
where we dropped the $B^2 \Lnorm{A^\fr{1}{2} u^\eta}{2}^2$ term from the left-hand side, split the $\nu \Lnorm{A u^\eta}{2}^2$ term, and applied Poincar\'e inequality, in order to remove the dependence on $B^2$ from the decay term.
Finally, we employ the Gr\"onwall Lemma~\ref{gronwall_lemma} with $c_1 = \min \left( C_6 \nu, 4C_4 \right)$, $c_2 = \min \left( \nu/2, C_5 \right)$, $G = H = 0$,
\begin{equation}
    F(t) = C \left( \Lnorm{u^\eta}{2}^2 \Lnorm{A^\fr{1}{2} u^\eta}{2}^2 + (1 + B^2) \Lnorm{\l^\fr{1}{2} q^\eta}{2}^2 \right), \quad  I(t) = C \Lnorm{q^\eta}{2}^4 + CB^2 \Lnorm{q^\eta}{2}^2,
\end{equation}
and $a=b=C_7\min \left( 1, \fr{1}{B^2 + \nu} \right)$. $F$ inherits the time integrability of $\Lnorm{A^\fr{1}{2} u^\eta}{2}^2$ and $\Lnorm{\l^\fr{1}{2} q^\eta}{2}^2$ from \eqref{decay_of_u_and_lambda_minus_half_q_and_time_integrability_of_A_half_u} and \eqref{decay_of_q_and_lambda_half_q_and_time_integrability} respectively, and $I$ inherits the exponential decay from \eqref{decay_of_q_and_lambda_half_q_and_time_integrability}.
Thus, we have
\begin{equation}
\label{final_decay_and_time_integrability}
    \Lnorm{A^\fr{1}{2} u^\eta}{2}^2 + \Lnorm{q^\eta}{4}^4 \leq e^{M^2} \left(1+ \frac{1}{\min \left(1, \frac{1}{B^2+\nu} \right)}\right) \G_5 e^{- C_8\min \left( 1, \fr{1}{B^2 + \nu} \right) t}, 
\end{equation}
and 
\begin{equation}
    \TInt{\left(\Lnorm{A u^\eta}{2}^2 + \Lnorm{\l^\fr{1}{2} (q^\eta)^2}{2}^2\right)} \leq e^{M^2} \left(1   + \fr{1}{\min\left(1, \fr{1}{B^2 + \nu} \right)}\right)^2 \G_6.
\end{equation}

With the estimates derived in Steps 1--3, one can establish the existence of solutions with the desired decay property \eqref{decay-thm}. Finally, in the next step, we prove the uniqueness of solutions.

\medskip

\textbf{Step 4. Uniqueness.}
Next, for the proof of uniqueness, we consider two solutions $(q_1, u_1)$ and $(q_2,u_2)$ of the system \eqref{charge_density_eq}-\eqref{div_free_eq}, with respective pressures $p_1$ and $p_2$, and we define the differences $u \coloneqq u_2 - u_1$, $q \coloneqq q_2 - q_1$, and $p\coloneqq p_2 - p_1$. These differences obey the following equations 
\begin{subequations}
    \begin{align}
    \partial_t u + u\cdot \na u_1+ u_2 \cdot \na u + \na p + B^2 u - \nu \D u &+ q R q_1 + q_2Rq \notag \\
    &= - Bq_1 u^\perp - Bqu_2^\perp + B R^\perp q, \label{uniq_u} \\
    \partial_t q + u \cdot \na q_2 + u_1 \cdot \na q + \l q &= B \Div{u^\perp}, \label{uniq_q} \\
    \Div{u} &= 0. \label{uniq_div_free}
\end{align}
\end{subequations}
We multiply \eqref{uniq_u} by $u$ and \eqref{uniq_q} by $\l^{-1} q$, integrate spatially over $\Omega$, and sum. We make use of the cancellation identity
\be\begin{aligned}
    \OInt{(u \cdot \na q_2) \l^{-1} q} &+ \OInt{q_2 Rq \cdot u} = 0, 
\end{aligned}\ee
and the fact that the pressure term $\OInt{\na p \cdot u}$, the cubic term $\OInt{(u_2 \cdot \na u) \cdot u}$, and the $B\OInt{q_1 u^\perp \cdot u}$ term vanish due to the divergence-free property of $u$ and $u_2$ and integration by parts, leaving us with the following evolution equation for $\l^{-\fr{1}{2}} q$ and $u$ in $L^2$,
\be\begin{aligned}
    \fr{1}{2} &\dt \left( \Lnorm{\l^{-\fr{1}{2}} q}{2}^2 + \Lnorm{u}{2}^2 \right) + \Lnorm{q}{2}^2 + B^2 \Lnorm{u}{2}^2 + \nu \Lnorm{\na u}{2}^2 \\
    &= - \OInt{(u \cdot \na u_1) \cdot u} - B \OInt{q u_2^\perp \cdot u} + B \OInt{R^\perp q \cdot u} - \OInt{q R q_1 \cdot u} \\
    &\quad \quad - \OInt{(u_1 \cdot \na q) \l^{-1} q} + B \OInt{(\Div{u^\perp}) \l^{-1} q}.
\end{aligned}\ee
Crucially, the advective term $\OInt{(u_1 \cdot \na q) \l^{-1} q}$ is critical and requires the enhanced commutator estimate we derived in Theorem \ref{commutator_theorem_2} (with the choice of $s=1, \, \kappa = 1/6$), resulting in the bound
\be\begin{aligned}
    \left| \OInt{(u_1 \cdot \na q) \l^{-1} q}\right| &= \left|\OInt{\l^{-\fr{1}{2}}(u_1 \cdot \na q) \l^{-\fr{1}{2}} q}\right| \\
    &= \left| \OInt{\left[ \l^{-\fr{1}{2}} (u_1 \cdot \na q) - u_1 \cdot \na \l^{-\fr{1}{2}} q \right] \l^{-\fr{1}{2}} q}\right| \\
    &\leq \fr{1}{6} \Lnorm{q}{2}^2 + C \left(1 + \Lnorm{\na u_1}{2}^{2} \right) \Lnorm{\D u_1}{2}^2 \Lnorm{\l^{-\fr{1}{2}} q}{2}^2,
\end{aligned}\ee 
after subtracting the vanishing term $\OInt{(u_1 \cdot \na \l^{-\fr{1}{2}}q) \l^{-\fr{1}{2}}q}$.

The remaining terms are of lower order and thus can be estimated more easily.
The quadratic terms can be directly bounded through the H\"older and Young inequalities, giving
\begin{equation}
    B \left| \OInt{R^\perp q \cdot u} \right| \leq \fr{1}{12} \Lnorm{q}{2}^2 + C B^2 \Lnorm{u}{2}^2,
\end{equation}
having used that the Riesz operator is an isometry in $L^2$, and
\begin{equation}
    B \left| \OInt{(\Div{u^\perp}) \l^{-1} q} \right| \leq \fr{1}{12} \Lnorm{q}{2}^2 + C B^2 \Lnorm{u}{2}^2.
\end{equation}
The cubic terms are similarly treated by employing the H\"older, Young, and Ladyzhenskaya inequalities, giving rise to the estimates
\begin{equation}
    \left| \OInt{(u \cdot \na u_1) \cdot u} \right| \leq \Lnorm{u}{4}^2 \Lnorm{\na u_1}{2} \leq \Lnorm{u}{2} \Lnorm{\na u}{2} \Lnorm{\na u_1}{2} \leq \fr{\nu}{6} \Lnorm{\na u}{2}^2 + C \Lnorm{\na u_1}{2}^2 \Lnorm{u}{2}^2, 
\end{equation}
and
\begin{equation}
    B \left| \OInt{q u_2^\perp \cdot u} \right| \leq B \Lnorm{q}{2} \Lnorm{u_2}{4} \Lnorm{u}{2}^\fr{1}{2} \Lnorm{\na u}{2}^\fr{1}{2} \leq \fr{1}{12} \Lnorm{q}{2}^2 + \fr{\nu}{6} \Lnorm{\na u}{2}^2 + C B^4 \Lnorm{u_2}{4}^4 \Lnorm{u}{2}^2.
\end{equation}
The boundedness of the Riesz operator in $L^4$ further allows us to bound the remaining term
\begin{equation}
    \left| \OInt{q Rq_1 \cdot u} \right| \leq C \Lnorm{q}{2} \Lnorm{q_1}{4} \Lnorm{u}{2}^\fr{1}{2} \Lnorm{\na u}{2}^\fr{1}{2} \leq \fr{1}{12} \Lnorm{q}{2}^2 + \fr{\nu}{6} \Lnorm{\na u}{2}^2 + C \Lnorm{q_1}{4}^4 \Lnorm{u}{2}^2.
\end{equation}
Thus, we reach the differential inequality
\be\begin{aligned}
    \dt \left( \Lnorm{\l^{-\fr{1}{2}} q}{2}^2 + \Lnorm{u}{2}^2 \right) &\leq C(1 + \Lnorm{\na u_1}{2}^2) \Lnorm{\D u_1}{2}^2 \Lnorm{\l^{-\fr{1}{2}} q}{2}^2 \\
    &\quad \ + C\left(B^2 + \Lnorm{q_1}{4}^4 + B^4 \Lnorm{u_2}{4}^4 + \Lnorm{\na u_1}{2}^2 \right) \Lnorm{u}{2}^2.
\end{aligned}\ee
By the regularity of the solutions established in the previous steps, the coefficient is integrable in time. The standard Grönwall inequality then gives the desired uniqueness property.

\end{proof}

\appendix
\section{Proof of Proposition \ref{heat_kernel_estimate}: Heat Kernel Estimate}
\label{kernel_appendix}
To prove Proposition \ref{heat_kernel_estimate}, we first consider the case when $|x-y| \geq \fr{d(x)}{10}$. We get, from the gradient bound \eqref{gradient_bound}, that the following bound holds for all $t > 0$,
\begin{equation}
    |(\na_x + \na_y) H_D| \leq C t^{-\fr{1}{2}-\fr{d}{2}} e^{-\fr{d(x)^2}{200Kt}} e^{-\fr{|x-y|^2}{2Kt}}. \label{gradient_bound_2}
\end{equation}
Hence, in view of \eqref{standard_integral_estimate}, \eqref{gradient_bound_2}, and a standard Gaussian integral on $\R^d$, we have
\be\begin{aligned}
\label{bound_1}
    &\int_0^\infty t^{-1+\fr{s}{4}} \left( \int_{|x-y|\geq \fr{d(x)}{10}} |x-y|^{q \alpha} \, |(\na_x + \na_y) H_D|^q \ dy \right)^\fr{1}{q} \ dt \\
    &\qquad \lesssim \int_0^\infty t^{-\fr{3}{2}+\fr{s}{4} - \fr{d}{2}} \, e^{-\fr{d(x)^2}{200Kt}} \left( \int_{|x-y|\geq \fr{d(x)}{10}} |x-y|^{q \alpha} \, e^{-\fr{q|x-y|^2}{2Kt}} \ dy \right)^\fr{1}{q} \ dt \\
    &\qquad \lesssim \int_0^\infty t^{-\fr{3}{2}+\fr{s}{4}-\fr{d}{2}} \, e^{-\fr{d(x)^2}{200Kt}} \left( \int_{\R^d} |x-y|^{q \alpha} \, e^{-\fr{q|x-y|^2}{2Kt}} \ dy \right)^\fr{1}{q} \ dt \\
    &\qquad \lesssim \int_0^\infty t^{-\fr{3}{2}+\fr{s}{4} - \fr{d}{2} + \fr{d}{2q} + \fr{\alpha}{2}} \, e^{-\fr{d(x)^2}{200Kt}} \ dt \\
    &\qquad \lesssim d(x)^{\alpha + \fr{s}{2} - d + \fr{d}{q} - 1},
\end{aligned}\ee
provided that $\alpha < 1 + d - \fr{d}{q} - \fr{s}{2}$ and $q\alpha>-d$. \\
Otherwise, if $|x-y| \leq \fr{d(x)}{10}$, it follows from Appendix 1 of \cite{constantin2016critical} that, for $t \leq d(x)^2$,
\begin{equation}
    |(\na_x + \na_y) H_D| \leq C t^{-\fr{1}{2} - \fr{d}{2}} e^{-\fr{d(x)^2}{Kt}}. \label{gradient_estimate_3}
\end{equation}
Making use of \eqref{standard_integral_estimate}, \eqref{gradient_estimate_3}, and integration via polar coordinates, we obtain that
\be\begin{aligned}
\label{bound_2}
    &\int_0^{d(x)^2} t^{-1 + \fr{s}{4}} \left( \int_{|x-y| \leq \fr{d(x)}{10}} |x-y|^{q\alpha} \, |(\na_x + \na_y) H_D|^q \ dy \right)^\fr{1}{q} \ dt \\
    &\qquad \lesssim \int_{0}^{d(x)^2} t^{-\fr{3}{2} + \fr{s}{4} - \fr{d}{2}} e^{-\fr{d(x)^2}{Kt}} \left(\int_{|x-y| \leq \fr{d(x)}{10}} |x-y|^{q\alpha} \ dy \right)^\fr{1}{q} \ dt \\
    &\qquad \lesssim \int_{0}^{d(x)^2} t^{-\fr{3}{2} + \fr{s}{4} - \fr{d}{2}} e^{-\fr{d(x)^2}{Kt}} d(x)^{\alpha + \fr{d}{q}} \ dt 
    = C \int_{0}^{d(x)^2} t^{-\fr{3}{2} + \fr{s}{4} - \fr{d}{2}} d(x)^{\alpha + \fr{d}{q}} e^{-\fr{d(x)^2}{K_1 t}} e^{-\fr{d(x)^2}{K_2 t}} \ dt \\
    &\qquad \lesssim \int_0^{d(x)^2} t^{-\fr{3}{2} + \fr{s}{4} - \fr{d}{2} + \fr{d}{2q} + \fr{\alpha}{2}} e^{- \fr{d(x)^2}{K_2 t}}\ dt \\
    &\qquad \lesssim d(x)^{\alpha + \fr{s}{2} - d + \fr{d}{q} - 1},
\end{aligned}\ee
having used the fact that $q\alpha > -d$ and $d(x)^{\alpha + \fr{d}{q}} e^{-\fr{d(x)^2}{K_1 t}} \leq C t^{\fr{\alpha}{2} + \fr{d}{2q}}$. \\
What is left to estimate is the integral
\begin{equation}
    \text{I} \coloneqq \int_{d(x)^2}^\infty t^{-1+\fr{s}{4}} \left( \int_{|x-y|\leq \fr{d(x)}{10}} |x-y|^{q \alpha} \, |(\na_x + \na_y) H_D|^q \ dy \right)^\fr{1}{q} \ dt.
\end{equation}
Using the gradient bound \eqref{gradient_bound}, we obtain
\begin{equation}
    |(\na_x + \na_y) H_D| \leq C t^{-\fr{1}{2} - \fr{d}{2}} e^{- \fr{|x-y|^2}{Kt}},
\end{equation}
which, in turn, allows us to estimate $\text{I}$, using the Gaussian integral, as
\be\begin{aligned}
\label{bound_3}
    \text{I} &\lesssim \int_{d(x)^2}^\infty t^{-\fr{3}{2} +\fr{s}{4} - \fr{d}{2}} \left( \int_{|x-y|\leq \fr{d(x)}{10}} |x-y|^{q\alpha} \, e^{-\fr{q|x-y|^2}{Kt}} \ dy \right)^\fr{1}{q} \ dt \\
    &\lesssim \int_{d(x)^2}^\infty t^{-\fr{3}{2} + \fr{s}{4} - \fr{d}{2} + \fr{d}{2q} + \fr{\alpha}{2}} \ dt \\
    &\lesssim d(x)^{\alpha + \fr{s}{2} -  d + \fr{d}{q} - 1},
\end{aligned}\ee
where we have used that $q\alpha>-d$.
Combining the integral bounds \eqref{bound_1}, \eqref{bound_2}, and \eqref{bound_3}, we reach the desired result \eqref{heat_kernel_result}.

\section{Gr\"onwall-Type Inequality}
\label{gronwall_section}
In this appendix, we prove the following Gr\"onwall-type inequality.
\begin{lemma} 
\label{gronwall_lemma}
Suppose $y$ obeys the differential inequality
\begin{equation}
\label{given_inequality}
    \dt y + c_1 y + c_2 z \leq F(t) y + G(t) H(t) + I(t),
\end{equation}
with $c_1> 0$ and $c_2>0$, where, for all $t \geq 0$,
\begin{equation}
\label{gronwall_conditions}
    y, z, F, G, H, I \geq 0, \ \TInt{F(s)} \leq A_1, \ \TInt{G(s)} \leq B, \ H(t) \leq A_2 e^{-at}, \ I(t) \leq A_3 e^{-bt},
\end{equation}
for some $A_1, A_2, A_3, B \geq 0$ and $a, b>0$. 
Let $\gamma=\min \left(c_1, a/2, b/2 \right)$. Then
\begin{equation}
\label{decay_bound}
    y(t) \leq e^{A_1} \left( y(0) + A_2 B + \fr{ A_3}{b - \gamma} \right) e^{-\gamma t} \coloneqq \G e^{- \gamma t},
\end{equation}
and
\begin{equation}
    \label{time_integrability_of_z}    
    \TInt{z(s)} \leq \fr{1}{c_2} \left( y(0) + A_1 \G + A_2 B + \fr{A_3}{b} \right).
\end{equation}
\end{lemma}
\begin{proof}
We define $\gamma = \min(c_1, a/2, b/2)$. In \eqref{given_inequality}, we can bound $c_1 y \geq \gamma y$ to get the differential inequality
\begin{equation}
\label{next_diff_inequality}
    \fr{dy}{dt} + (\gamma - F(t)) y + c_2 z \leq G(t) H(t) + I(t).
\end{equation}
We choose the integrating factor to be
\begin{equation}
    \mu(t) = \exp \left( \TInt{(\gamma - F(s))}\right),
\end{equation}
which satisfies the following two inequalities
\begin{equation}
\label{integrating_factor_properties}
    \mu(t) \leq e^{\gamma t}, \qquad \mu(t) \geq e^{\gamma t-A_1},
\end{equation}
due to the non-negativity of $F$ (giving that $\exp \left( - \TInt{F(s)} \right) \leq 1$) and the integral boundedness condition on $F$ from \eqref{gronwall_conditions}.
We multiply \eqref{next_diff_inequality} by the integrating factor to get
\begin{equation}
    \dt (\mu y) \leq \mu(t) G(t) H (t) + \mu(t) I(t),
\end{equation}
having dropped the $z$ term on the left-hand side. Integrating in time from $0$ to $t$, we have
\be\begin{aligned}
\label{intermediary_gronwall}
    y(t) \mu(t) &\leq y(0) + \TInt{\mu(s) G(s) H(s)} + \TInt{\mu(s) I(s)} \\
    &\leq y(0) + A_2 \TInt{e^{(\gamma-a) s} \, G(s)} + A_3 \TInt{e^{(\gamma -b)s}} \\
    &\leq y(0) + A_2 \TInt{G(s)} + \fr{A_3}{\gamma - b} \left( e^{(\gamma - b) t} - 1 \right) \\
    &\leq y(0) + A_2 B + \fr{A_3}{b - \gamma},
\end{aligned}\ee
having used the imposed conditions \eqref{gronwall_conditions}, the upper bound on the integrating factor from \eqref{integrating_factor_properties}, and that $e^{(\gamma - a) s}, \,e^{(\gamma - b)s} \leq 1$. Hence, we deduce the decaying bound in \eqref{decay_bound} by dividing by the integrating factor and using the fact that $1/\mu(t) \leq e^{A_1 - \gamma t}$ (by \eqref{integrating_factor_properties}). 

Finally, what is left to show is the boundedness of the time integral of $z$. 
We integrate \eqref{given_inequality} in time from $0$ to $t$ and drop the $y(t)$ and $c_1\TInt{y(s)}$ terms on the left-hand side to obtain
\be\begin{aligned}
    c_2 \TInt{z(s)} &\leq y(0) + \TInt{F(s) y(s)} + \TInt{G(s) H(s)} + \TInt{I(s)} \\
    &\leq y(0) + \G \TInt{F(s) e^{-\gamma s}} + A_2 \TInt{G(s) e^{-as}} + A_3 \TInt{e^{-bs} }\\
    &\leq y(0) + \G\TInt{F(s)} + A_2 \TInt{G(s)} + A_3 \int_0^\infty e^{-bs} \, ds \\
    &\leq y(0) + A_1 \G + A_2 B + \fr{A_3}{b},
\end{aligned}\ee
where we employed the derived decaying bound on $y(t)$ and the given conditions \eqref{gronwall_conditions}. This completes the proof of \eqref{time_integrability_of_z}.
\end{proof}

\vspace{0.5cm}

{\bf{Acknowledgments.}} E.A. and J.G. were partially supported by the University Research Board (URB) of the American University of Beirut under Grant No. 104752.  Q.L. was partially supported by the Simons Foundation grant (SFI-MPS-TSM-00013384).

\vspace{0.5cm}

{\bf{Data Availability Statement.}} The research does not have any associated data.

\vspace{0.5cm}

{\bf{Conflict of Interest.}} The authors declare that they have no conflict of interest.

\bibliographystyle{plain}
\bibliography{bibliography}

\end{document}